\documentclass[12pt]{article}

\usepackage{graphicx}
\usepackage{mathptmx}
\usepackage{amsmath,amssymb,amsfonts}

\newcommand{\oM}{\overline{M}}
\newtheorem{theo}{Theorem}

\newtheorem{corol}{Corollary}

\newcommand{\sm}[1]{\mbox{\small $#1$}}

\newcommand{\La}[1]{\mbox{\Large $#1$}}
\newcommand{\LA}[1]{\mbox{\LARGE $#1$}}

\newcommand{\h}{\mbox{\lie h}}                  
\newcommand{\R}{\mathbb{R}}          

\newfont{\lie}{eufm10 at 12pt}
\newfont{\field}{msbm10 at 11pt}

\begin{document}
\title{\bf \small  MEAN CURVATURE FLOW 
AND BERNSTEIN-CALABI RESULTS FOR SPACELIKE GRAPHS}
\author{Guanghan Li$^{1,2}$ and Isabel M.C.\ Salavessa$^{2,\dag}$}
\date{}
\maketitle ~~~\\[-10mm]
{
\footnotesize $1$ School of Mathematics and Computer Science,
Hubei University, Wuhan, 430062, P. R. China,\\[-1mm]
e-mail: liguanghan@163.com}\\[1mm]
{\footnotesize $2$ Centro de F\'{\i}sica das Interac\c{c}\~{o}es
Fundamentais, Instituto Superior T\'{e}cnico, Technical University\\[-1mm]
of Lisbon, Edif\'{\i}cio Ci\^{e}ncia, Piso 3, Av.\ Rovisco Pais,
1049-001 Lisboa, Portugal;\\[-1mm]
$\dag$e-mail: isabel.salavessa@ist.utl.pt}
\markright{\sl\hfill Li -- Salavessa \hfill}\\[3mm]
{\em  \footnotesize
Contribution to the  VIII International Colloquium on Differential Geometry,
Santiago de Compostela, Spain, 7-11 July 2008,
 satellite event of the 5th European Congress of Mathematics.}
\begin{abstract} This is  a survey of our work 
  on spacelike
graphic submanifolds in pseudo-Riemannian products, namely on 
Heinz-Chern and 
Bernstein-Calabi
results and on the mean curvature flow, with applications to the
homotopy of maps between Riemannian manifolds.
\end{abstract}
\section{Introduction}
\renewcommand{\thesection}{\arabic{section}}
\renewcommand{\theequation}{\arabic{equation}}
\setcounter{equation}{0}
It has been an important issue in geometry and in topology
to determine when a map $f:\Sigma_1\to \Sigma_2$ between manifolds
can be  homotopically deformed to a constant one.
If each $\Sigma_i$  has a Riemannian structure $g_i$,
the curvature of these spaces may give an answer.
This is particularly more complex if $\Sigma_i$ are both compact.
For $\Sigma_i$ noncompact, by a famous result due to Gromov (\cite{grom}),
$\Sigma_i$ admits a Riemannian metric of negative sectional curvature
and also one of positive sectional curvature. In each of these cases,
if $\Sigma_i$ is  complete and simply connected,
 then  $\Sigma_i$ is  diffeomorphic to a  contractible space, 
by the Cartan-Hadmard  theorem and by a result of Cheeger and Gromoll, 
respectively (see in \cite{cheebin}).
If this is the case for one of the $\Sigma_i$
then  $f$ is obviously  homotopically trivial.

A deformation problem of an initial map can be handled using some
geometric evolution equation, obtaining homotopic deformations 
  of a certain type and with
geometrical and analytical meaning, namely, giving at infinite time
a solution of a certain partial differential equation.
We recall the great discovery of Eells and Sampson (\cite{es}),
a first example of this kind,
on using the heat flow to deform a map
to an harmonic one:
\begin{theo}[Eells and Sampson (1964)]
 If $\Sigma_1$ and $\Sigma_2$ are  closed and 
$\Sigma_2$ has nonpositive sectional curvature then
$f$ is homotopic to
a harmonic map $f_{\infty}$. Furthermore, if the Ricci tensor  of $\Sigma_1$ 
is nonnegative
then $f_{\infty}$ is totally geodesic and if it is positive somewhere, then
$f_{\infty}$ is constant. 
\end{theo}
\noindent
The last part of this theorem can be seen as a Bernstein-type theorem, and
it was obtained from a Weitzenb\"{o}ck formula for
the Laplacian of $\|df_{\infty}\|^2$. We recall that Bernstein-type theorems
are theorems that give conditions that
ensures that a solution of  certain P.D.Es.\ with geometrical meaning,
 must be a   a "trivial" solution, 
as for example
a totally geodesic or a constant map.

In this note, a survey of our main results in \cite{lisa1,lisa2,sal},
we will show how to use the mean curvature flow and a  Bernstein-Calabi type
result
for spacelike graphs to
obtain a  deformation of a map between Riemannian manifolds to a 
totally geodesic  or  a constant one.

The Bernstein-Calabi result is obtained by computing the Laplacian of
a positive geometric quantity, the hyperbolic cosine of the
hyperbolic angle of a spacelike graph, and analyzing the sign
of this Laplacian,
 based on an idea of  Chern \cite{chern} of computing
a similar quantity in the Riemannian case.

Furthermore, we  also will show that
under somehow more general curvature conditions as in the above theorem, we can
obtain a direct proof of the homotopy to a constant map, with no need to
use a Bernstein-type result.
This approach was started by
Wang \cite{w1} for the graph $\Gamma_f$ of the map $f$, considered
as a submanifold of the Riemannian product $\Sigma_1\times \Sigma_2$
of closed spaces with constant sectional curvature, and take its
mean curvature flow and show that under certain conditions
the flow preserves the graphic structure of the submanifold
and converges to the graph of a constant map.
The main difference with our approach is that we
 consider the pseudo-Riemannian structure on $\Sigma_1\times \Sigma_2$
instead the Riemannian one. Our assumption on $\Gamma_f$ to be
a  spacelike submanifold is essentially identical to the
assumptions on the eigenvalues of $f^*g_2$  imposed in \cite{w1,tw}  
in the corresponding Riemannian setting. 
Our advantage is that the pseudo-Riemannian setting
turns out to be a more natural one, since it allows less restrictive
assumptions on the curvature tensors (and that include the case
of any negative sectional curvature for $\Sigma_2$) 
and on the map $f$ itself after a suitable rescaling of the
metric of $\Sigma_2$,
and long time existence and convergence of
the flow are easier  obtained. 
In \cite{w1} it is necessary to use a White's regularity theorem,
based on a monocity formula due to Huisken, to detect possible singularities 
of the mean curvature flow, while in the pseudo-Riemannian case, because
of good signature in the evolution equations, we have better regularity.
This will become clear in equations (1), (3)
and (4) below.\\

Let $(\Sigma_1,g_1)$ and $(\Sigma_2,g_2)$ be Riemannian manifolds of dimension
$m\geq 2$ and $n\geq 1$ respectively, and of sectional curvatures $K_i$ 
and Ricci tensors
$Ricci_i$. On $\oM=\Sigma_1\times \Sigma _2$ we
consider the  pseudo-Riemannian metric $\overline{g}=g_1-g_2$. We assume
$\Sigma_1$ oriented.
Given a map $f$, we assume  the graph, $\Gamma_f:\Sigma_1\to \oM$,
 $\Gamma_f(p)=(p,f(p))$, is a spacelike submanifold
that is, $g:=\Gamma_f^*\bar{g}=g_1-f^*g_2$ is a Riemannian 
metric on $\Sigma_1$. Thus, the eigenvalues
of $f^*g_2$, at $p\in M$,
$\lambda_1^2\geq \ldots\geq \lambda_m^2\geq 0$, are bounded from above
by $1-\delta(p)$, where $0< \delta(p)\leq 1$ is
a constant depending on $p$. The \em hyperbolic angle \em
$\theta$ of $\Gamma_f$ is given by one of the equivalent definitions:
$$ \cosh\theta=\frac{1}{\sqrt{\Pi_i(1-\lambda_i^2)}}=
Vol_{\Sigma_1}(\pi_1(e_1),\ldots,\pi_1(e_m))=\frac{Vol_{(\Sigma_1, g_1)}}{
Vol_{(\Sigma_1,g)}}$$
where $\pi_1:\overline{M}\to \Sigma_1$ is the projection and
$e_i$ is a direct o.n. basis of $T\Gamma_f$,
and $Vol_{(\Sigma_1,g)}$ is the volume element of
$(\Sigma_1,g)$. Then $\cosh\theta\equiv 1$
iff $f$ is constant, that is $\Gamma_f$ is a slice.
\section{Bernstein-Calabi and Heinz-Chern type results}
The classic Bernstein theorem says that an entire minimal graph in
$\R^3$ is a plane. This result was generalized to codimension one graphs
in $\R^{m+1}$ for
$m\le 7$, and for higher dimensions and codimensions under additional
 conditions by many other authors.  
Calabi (\cite{cal}) considered  the same problem for
the maximal (the mean curvature
$H=0$) spacelike hypersurfaces $M$ in the Lorentz-Minkowski space
$\R_1^{m+1}$ with the
metric
$ds^2=\sum _{i=1}^m(dx_i)^2-(dx_{m+1})^2.$
If $M$ is given by the graph of a function $f$ on $\R^m$
with $|Df|<1$, the equation  $H=0$ has the form
$$\sum _{i=1}^m\frac {\partial }{\partial x_i}
\left(\frac {{\partial f}/{\partial
x_i}}{\sqrt{1-|Df|^2}}\right)=0.$$
Calabi
showed that for $ m\le 4$, the graph of any entire solution to the
above equation is a hyperplane. The same conclusion was established
  by Cheng and Yau (\cite{cy})
for any $m$.
A further generalization of this problem to
$\mathbb{R}^{m+n}_n$ has been obtained by some authors
(see for instance in \cite{jx}).

 Another natural generalization is to 
consider  maximal 
spacelike graphic submanifolds in a non flat ambient
space and in higher codimension. We consider a spacelike graph
$\Gamma_f$, for a map $f:\Sigma_1\to \Sigma_2$.

We can take
$a_i$ an orthonormal basis of $T_p\Sigma_1$ and $e_{\alpha}$ of
$T_{f(p)}\Sigma_2$, $1\leq i\leq m$,  $ m+1\leq \alpha\leq m+n$,
such that $df(a_i)=-\lambda_{i}a_{m+i}$ ($\lambda_i=0$ if $i>n$). Then
$e_i=(1-\lambda_i^2)^{-1/2}(a_i+\lambda_ia_{m+i})$ 
and $e_{m+i}=(1-\lambda_i^2)^{-1/2}(a_{m+i}+\lambda_ia_i)$,
$e_{\alpha}=a_{\alpha}$ if $\alpha>2m$,
define
o.n.bs. \ of $T_{(p,f(p))}\Gamma_f$, and of
the normal bundle at $(p,f(p))$ respectively. Assuming $\Gamma_f$ has parallel
mean curvature, in this basis we have
\begin{eqnarray}
\lefteqn{\Delta \cosh\theta=
 \cosh\theta \LA{\{} ||B||^2+2\sum_k\sum _{i<j}\lambda
_i\lambda _jh_{ik}^{m+i}h_{jk}^{m+j}-2\sum_k\sum _{i<j}\lambda
_i\lambda _jh_{ik}^{m+j}h_{jk}^{m+i}} \nonumber\\
&&+\sum_{i}\La{(}{\sm{\frac{ \lambda _i^2}{(1-\lambda
_i^2)}}}Ricci_1(a_i,a_i)+\sum_{j\neq i}\sm{\frac{ \lambda _i^2\lambda_j^2}
{(1-\lambda_i^2)(1-\lambda _j^2)}} [K_1(P_{ij})-K_2(P'_{ij})]\La{)}
\LA{\}}
\end{eqnarray}
where $P_{ij}=span\{a_i,a_j\}$ and 
 $P'_{ij}=span\{a_{m+i},a_{m+j}\}$.
Here $h^{\alpha}_{ij}$ are the components of the
second fundamental form  $B$ of $\Gamma_f$ in the basis $e_i$, $e_{\alpha}$.
 
\begin{theo}[{\cite{sal,lisa1}}]
Let $M=\Gamma_f$ be a spacelike graph submanifold of $\overline{M}$
with parallel mean curvature vector.
We assume for each $p\in \Sigma_1$, $Ricci_1(p)\geq 0$ and
for any two-dimensional planes  $P\subset T_p\Sigma_1$, $P'\subset
T_{f(p)}\Sigma_2$, 
  $K_1(P)\ge K_2(P')$. We have:\\[2mm]
$(i)$~If $n=1$ and $\cosh\theta\leq o(r)$ when $r\to +\infty$, where
$r$ is the distance function to a point  $p\in (\Sigma_1,g_1)$, 
and $\Sigma_1$ is complete, then
$M$ is maximal.\\[1mm]
$(ii)$~  If $M$ is compact, then it is totally geodesic. Moreover, if
$Ricci_1(p)>0$ at some point, then $M$
is a slice, that is $f$ is constant;\\[1mm]
$(iii)$~  If $M$ is complete, noncompact, and $K_1$,
$K_2$ and  $\cosh\theta$ are bounded, then $M$ is maximal. \\[1mm]
$(iv)$~ If $M$ is a complete maximal spacelike surface, 
then $M$ is totally geodesic. Moreover, $(a)$
if $K_1(p)>0$ at some point $p\in M$,  then $M$ is a slice;
$(b)$ If $\Sigma_1=\mathbb{R}^2$ and $\Sigma_2=\mathbb{R}^n$, 
then $M$ is a plane; $(c)$ if $\Sigma_1$ is flat
and $K_2<0$ at some point $f(p)$, then either
$M$ is a slice or the image of $f$ is a geodesic of $\Sigma_2$.
\end{theo} 
\noindent
We obtain $(i)$  by applying  a Heinz-Chern inequality derived  in
\cite{sal}, for the absolute norm of $H$
$$m\|H\|\leq \sup_{B_r(p)}\cosh\theta\,\, \h(B_r(p)),$$
where $\h(B_r(p))=\inf_D V_{m-1}(\partial D)/V_m(D)$,
 is the Cheeger constant of the open geodesic ball
of center $p$ and  radius $r$, where $D$ runs all over the bounded domains
of the ball with smooth boundary $\partial D$.
Since $Ricci_1\geq 0$,
$\h(B_r(p))\leq C/r$, when $r\to +\infty$,
where $C>0$ is a constant. For $\Sigma_1$
 the $m$-hyperbolic space (with non-zero Cheeger constant),
we give examples in
\cite{sal}
of  foliations of $\mathbb{H}^m\times\mathbb{R}$ by
complete spacelike graphic hypersurfaces
with bounded hyperbolic angle and with constant mean curvature any real $c$,
the same for all leaves, or parameterized by the leaf.

The proof of $(ii)$ and  $(iii)$ consists on showing that, 
under the curvature conditions, one has
$\Delta \ln\cosh\theta$ $ \geq \delta\|B\|^2$, where $\delta>0$ is a 
constant that
does not depend on $p$  and in $(iii)$ showing
 that the Ricci tensor of $M$ is bounded from below, and applying
 the Omori-Cheng-Yau
maximum principle for noncompact manifolds.
$(i)$ and $(iii)$ are obtained
by different approaches.
If $M$ is a maximal Riemannian surface, 
$(iv)$ gives 
 a  generalization of  the
Bernstein type theorem of Albujer-Al\'{\i}as \cite {aa} for maximal
graphic spacelike surfaces in a Lorentzian three manifold  $\Sigma_1\times
\mathbb{R}$ to 
higer codimension. As in \cite{chern, aa} the proof is based
on a parabolicity argument for surfaces with nonnegative Gauss curvature.
In fact, in this case,  we have that
$\Delta\left(\frac{1}{\cosh\theta}\right)\leq 0$
and  the Gauss curvature of $M$ satisfies
$$K_M=\sm{\frac{1}{(1-\lambda_1^2)(1-\lambda_2^2)}}[K_1-\lambda_1^2\lambda_2^2
K_2(a_3,a_4)]+\sum_{\alpha}[(h_{11}^{\alpha})^2+(h_{12}^{\alpha})^2]
\geq 0.$$
The conclusion that $B=0$ comes from analyzing the vanishing of the term
involving the components of $B$ in the expression of $\Delta(1/\cosh\theta)$.
Our proof for $(iv)(b)$, gives a simpler proof of the same result of
Jost and Xin \cite{jx} for de case of surfaces, but using their result
that any entire maximal
graph in $\mathbb{R}^{m+n}_n$ is complete.\\

We also derive in \cite{lisa2}
a Simons' type identity for the absolute norm of the second 
fundamental form $\|B\|^2$ of a spacelike submanifold $M$ of any
pseudo-Riemannian manifold $\oM$, 
\begin{eqnarray*}
\Delta||B||^2&=&2||\nabla B||^2+\sm{\sum}_{ij\alpha}2
h_{ij}^{\alpha}H_{,ij}^{\alpha}-\sm{\sum}_{ij\alpha}2h_{ij}^{\alpha}
[\sm{\sum}_k(\bar{\nabla}_j\bar{R})_{kik}^{\alpha}+\sm{\sum}_k(\bar{\nabla}
_k\bar{R})_{ijk}^{\alpha}]\nonumber\\
&&+\sm{\sum}_{ij\alpha\beta}2\{\sm{\sum}_k(4\bar{R}_{\beta
ki}^{\alpha}h_{kj}^{\beta}h_{ij}^{\alpha}-\bar{R}_{k\beta
k}^{\alpha}h_{ij}^{\alpha}h_{ij}^{\beta})+\bar{R}_{i\beta
j}^{\alpha}H^{\beta}h_{ij}^{\alpha}\}\nonumber\\
&&+\sm{\sum}_{ijkl\alpha}4(\bar{R}_{ijk}^lh_{ij}^{\alpha}h_{kl}^{\alpha}
+\bar{R}_{kik}^lh_{lj}^{\alpha}h_{ij}^{\alpha})
-\sm{\sum}_{ijk\alpha\beta}2h_{ij}
^{\alpha}h_{jk} ^{\alpha}h_{ki} ^{\beta}H^{\beta}\nonumber\\
&&+2\sm{\sum}_{ij\alpha\beta}\La{(}\sm{\sum}_{k}(h_{ik} ^{\alpha}h_{jk}
^{\beta}-h_{ik} ^{\beta}h_{jk} ^{\alpha})\La{)}^2+2\sm{\sum}_{\alpha
\beta}(\sm{\sum}_{ij}h_{ij}^{\alpha}h_{ij}^{\beta})^2,
\end{eqnarray*}
\section{The mean curvature flow}
The mean curvature flow of an immersion $F_0:M\to \oM$ 
is a family of immersions $F_t:M\to \oM$  defined
in a maximal interval $t\in [0,T)$ evolving according to
\begin{equation}
\left\{\begin{array}{l}
\frac{d}{dt}F(x,t)=H(x,t)=\Delta_{g_t}F_t(x)\\
F(\cdot,0)=F_0
\end{array}\right.
\end{equation}
where $H_t$ is the mean curvature of $M_t=F_t(M)=(M,g_t=F^*_t
\bar{g})$. 
The mean curvature flow of hypersurfaces 
in a Riemannian manifold has been extensively studied.
Recently, mean curvature flow of submanifolds with
higher co-dimensions has been paid more attention.
In \cite{w1}, the graph mean curvature
flow is studied in Riemannian product manifolds, and it is proved long-time
existence and convergence of the  flow under
suitable conditions.
When $\oM$ is a pseudo-Riemannian manifold, it is considered  the
mean curvature flow of spacelike submanifolds. This flow
for spacelike hypersurfaces has also been largely studied, but very little
is known on mean curvature flow in higher codimensions except
in a flat space  $\mathbb{R}^{n+m}_n$ \cite{xin}.
In \cite{lisa2} we consider (2) with $\oM$ any pseudo-Riemannian
manifold and $F_0$ any spacelike submanifold, and 
 we derive
 the evolution of the following quantities at a given  point $(x,t)$
with respect to 
an o.n. frame
$e_{\alpha}$ of of the normal bundle of $M_t$
and a coordinate chart $x_i$ of $M$,
normal at  $x$ relatively to the  metric $g_t$ 
{ $$\begin{array}{ccl}
 \frac{d}{dt}g_{ij} &=& 2 H^{\alpha} h^{\alpha}_{ij}\\[3mm]
\frac{d}{dt}Vol_{M_t} &=&\|H\|^2Vol_{M_t}\\[3mm]
\frac d{dt}||B||^2
&=&\Delta ||B||^2-2||\nabla B||^2+\sm{\sum}_{ij\alpha}2h_{ij}^{\alpha}\La{(}
\sm{\sum}_k(\bar{\nabla}
_j\bar{R})_{kik}^{\alpha}+(\bar{\nabla}
_k\bar{R})_{ijk}^{\alpha}\La{)}\\[2mm]
&&-{\sm{\sum}}_{ijk\alpha\beta}2(4\bar{R}_{\beta
ki}^{\alpha}h_{kj}^{\beta}h_{ij}^{\alpha}-\bar{R}_{k\beta
k}^{\alpha}h_{ij}^{\alpha}h_{ij}^{\beta})
+{\sm{\sum}}_{ijkl\alpha}4(\bar{R}_{ijk}^lh_{ij}^{\alpha}h_{kl}^{\alpha}
+\bar{R}_{kik}^lh_{lj}^{\alpha}h_{ij}^{\alpha})\\[2mm]
&&-2\sm{\sum} _{ij\alpha\beta}\La{(}\sm{\sum}_{k}(h_{ik} ^{\alpha}h_{jk}
^{\beta}-h_{ik} ^{\beta}h_{jk} ^{\alpha})\La{)}^2-2\sm{\sum}_{\alpha,
\beta}(\sm{\sum}_{ij}h_{ij}^{\alpha}h_{ij}^{\beta})^2.
\end{array} $$}

In this section we assume $(\Sigma_1,g_1)$ closed and  $(\Sigma_2,
g_2)$ complete, and  the curvature
tensor $R_2$ of $\Sigma_2$ and all its covariant derivatives are bounded.
We consider the mean curvature flow on the pseudo-Riemannian manifold
$\oM=\Sigma_1\times \Sigma_2$, when the initial immersed submanifold
$F_0=\Gamma_f:M=\Sigma_1\to \oM$ is a spacelike graph submanifold.
Furthermore, we assume, as in theorem 2
$$Ricci_1(p)\geq 0,~~~\mbox{ and}~~
K_1(p)\ge K_2(q) ~~~~\forall p\in
\Sigma _1, q\in \Sigma _2.$$
This means that either $K_1(p)\geq K_2^+(q)
=\max\{K_2(q),0\}$, or $Ricci_1(p)\geq 0$ and 
$K_1(p)(P)$ $<0$ for some
two-plane $P$ and $K_1(p)\geq K_2(q)$, with $K_2(q)<0$,  $\forall p,q$.
\\

We recall the main steps of  \cite{lisa2}.

For $t>0$ sufficiently small,
 $F_t$ is near $F_0$ and so it is a spacelike graph with
$\lambda_i^2(t)\leq 1-\delta(t).$ We derive the 
evolution of the hyperbolic angle
\begin{eqnarray}
\lefteqn{\frac d{dt}\ln(\cosh\theta) =\Delta \ln(\cosh\theta)-
\La{\{}\underbrace{||B||^2-\sum _{k,i}\lambda _i^2(h_{ik}^{m+i})^2-2\sum
_{k,i<j}\lambda
_i\lambda _jh_{ik}^{m+j}h_{jk}^{m+i}}_{\begin{array}
{c}\geq \delta(t)\|B\|^2\\ ~~\end{array}}
\La{\}}}\nonumber \\[-2mm]
&&\!\!-\sum_i\lambda_i^2\La{(}\frac{1}{(1-\lambda
_i^2)}\underbrace{Ricci_1(e_i,e_i)}_{\geq 0}
+\sum_{i\neq j}\frac{\lambda_j^2}{(1-\lambda_i^2)(1-\lambda_j^2)}
\underbrace{[K_1(P_{ij})-K_2(P'_{ij})]}_{\geq 0}\La{)}.~~~
\end{eqnarray}
Therefore, $\frac d{dt}\ln(\cosh\theta)
\leq \Delta \ln(\cosh\theta)-\delta(t)\|B\|^2 \leq \Delta \ln(\cosh\theta)$,
and by the maximum principle for parabolic equations, 
$\max_{\Sigma_1}\cosh\theta_t$ is a 
nondecreasing function on $t$, 
and in particular  $F_t$ remains a spacelike graph
$F_t=\Gamma_{f_t}$ for a smooth map $f_t:\Sigma_1\to \oM$.
On what follows, $c_i$ denotes positive constants.
We may take a uniform bound $\delta=\delta(0)$, such that $\lambda_i^2(t)
\leq 1-\delta$ for all $t$ as long as the flow exists. Consequently
$g_t= g_1-f^*_tg_2$ are uniformly equivalent metrics on $\Sigma_1$ and
$Vol_{M_t}$  are uniformly bounded,
and from the above evolution equations 
$Vol_{M_t}=e^{\int_0^t\|H_s\|^2ds}Vol_{M_0}$,  
what implies $\int_0^T\sup_{\Sigma_1}\|H_t\|^2dt <c_0$.
From the evolution equations one gets
\begin{equation}
\frac d{dt}||B||^2~{\le}~
\Delta ||B||^2+c_1||B||+c_2||B||^2-\frac
2n||B||^4 ~{\le}~ \Delta ||B||^2
{-\frac 1n||B||^4}+c_3. 
\end{equation}
This is the point where regularity theory is better in the
pseudo-Riemannian setting than the Riemannian one
(note the negative coefficient of the highest power of $\|B\|$,
that holds in the pseudo-Riemannian case and not in the Riemannian
case).
From the above inequality
we may use  a result of Ecker and Huisken
\cite{eh} to conclude that $\|B\|^2$ is uniformly bounded. 
From this inequality we 
may apply an interpolation formula due to Hamilton and applying
parabolic maximum principles
we conclude $\|\nabla^kB\|^2$ is uniformly bounded for all $k$.\\
 
For each $t$ it is defined on $\oM$
a Riemannian metric  $\hat{g}_t=\bar{g}_{|_{T_pM_t}} - 
\bar{g}_{|_{T_pM_t^{\bot}}}$ that makes $e_i,e_{\alpha}$ an orthonormal basis.
These metrics defined along the flow are uniformly equivalent to the natural 
Riemannian metric $\bar{g}^+=g_1+g_2$ of $\oM=\Sigma_1\times \Sigma_2$, for
we have some positive constants $c(\delta)$ and $c'(\delta)$, 
depending only on $\delta$, such that $ c(\delta){\bar{g}_+}\leq
{\hat{g}}\leq c'(\delta){\bar{g}_+}$
holds. 
We observe that the Levi-Civita connections of $(\oM,\bar{g}_+)$
and of $(\oM,\bar{g})$ are
the same and   
$\|\bar{\nabla}B\|^2_{\hat{g}} \leq c_{22}\|B\|^4+
\|\nabla B\|^2.$ By induction on $k$ we see that $\bar{\nabla}^kB$ are
 $\hat{g}$ and so $\bar{g}_+$-uniformly bounded for
all $k\geq 0$, that is all derivatives of $B$ in $\oM$ are
also bounded for the Riemannian structure.
Then we can
apply Schauder theory, by embedding isometrically $(\Sigma_i,
g_i)$ into an Euclidean space $\mathbb{R}^{N_i}$. 
The spaces $C^{k+\sigma}(\Sigma_1,\oM)$, $k\in \mathbb{N}$,
$0\leq \sigma <1$ are Banach manifolds and can be seen as closed subsets
of the Banach space $C^{k+\sigma}(\Sigma_1,\mathbb{R}^{N_1+N_2})$ with
the H\"{o}lder norms.
Equation (2) in local
coordinates is of the form
$${\sum_{ij}a_{ij}\frac{\partial^2 F^{a}}{
\partial x_i\partial x_j} -\sum_k b_k\frac{\partial F^{a}}{
\partial x_k}= \bar{G}(x,t)^{a}+\frac{dF^{a}}{dt}}$$
where
 $a_{ij}=g^{ij}$, 
$b_k=g^{ij}\Gamma^l_{ij}$, $
\bar{G}(x,t)^{a}=(\bar{\Gamma}^a_{bc}\circ F_t)
\frac{\partial F^b}{\partial x_i}
\frac{\partial F^c}{\partial x_j}$.
From the uniform bounds of $\bar{\nabla}^kB$  and of
 $\bar{\nabla}^kH$ we have that
the coefficients $a_{ij},b_j$ are  $C^{k-1+\sigma}(\Sigma_1)$- uniformly
bounded, and if $F_t$ lies on a compact set of $\oM$ then 
$$ \|F(\cdot, t)\|_{C^{1+\sigma}(\Sigma_1,\oM)}
\leq c_{-1},~~~~~\|F(\cdot, t)\|_{C^{2+k+\sigma}(\Sigma_1,\oM)}\leq c_{k},
~~~k\geq 0$$
for some
positive constants $c_i$ that do not depend on $t$. 
Standard use of Ascoli-Arzela's theorem
to $F_t$ leads to the conclusion that $T=+\infty$ (by assuming $T<+\infty$
one has $F_t=F_0+\int_0^tH$ lies in a compact set and
gets an extension of the maximal solution $F_t$ to
$t=T$, what is a contradiction). We also note that the assumption
of $R_2$  and its derivatives to be bounded is necessary to guarantee
the existence of a maximal solution of the flow, as well
the trick of DeTurck can also be applied in the pseudo-Riemannian case
like in the Riemannian case, to reparametrize
$F_t$ in a suitable way to convert the above system in one of strictly
parabolic equations (see \cite{z} p. 17). This is necessary since the
coefficients $b_k$ also depend on the second derivatives of $F_t$, 
 and so it can give a degenerated system.
\begin{theo}[{\cite{lisa2}}] 
The mean curvature flow  of the spacelike
graph of $f$ remains a spacelike graph of a map
$f_t:\Sigma_1\rightarrow \Sigma_2$ and exists for all time $t\geq 0$.
\end{theo}
\noindent
Since $\int_0^{+\infty}\sup_{\Sigma_1}\|H_t\|^2 dt\leq c_{12}$, then
$\exists  t_N\rightarrow +\infty$  such that $H_{t_N}\to 0$. 
Assuming $f_t$ lies in a compact set
of $\Sigma_2$ we obtain a subsequence $F_{t_n}$ that $C^{\infty}$-converges
at infinity
to a map $F_{\infty}\in C^{\infty}(\Sigma_1,\oM)$, necessarily
a spacelike graph of a map $f_{\infty}\in 
C^{\infty}(\Sigma_1,\Sigma_2)$, and maximal, for $H_{\infty}=0$.
From Bernstein theorem 2, we conclude
\begin{theo}[{\cite{lisa2}}]  
If  $\Sigma_2$ is also compact there is a sequence $t_n\to +\infty$
such that the sequence $\Gamma_{f_{t_n}}$ of the
  flow  converges at infinity to a spacelike graph $\Gamma_{f_{\infty}}$
of a totally
geodesic map $f_{\infty}$, and if $Ricci_1(p)>0$ at some point $p\in \Sigma_1$, the sequence converges to a slice.
\end{theo}
Finally we consider the case $Ricci_1>0$ everywhere. In this case we 
will see that we can droop the compactness assumption of $\Sigma_2$.
From (3)
  $$\frac d{dt}\ln (\cosh\theta)\le \Delta \ln (\cosh\theta)
-c_{15}{\sum_i\lambda_i^2},$$
what implies
${
\frac d{dt}\ln (\cosh\theta)\le \Delta \ln (\cosh\theta)
-c_{15}\left(1-\frac
1{\cosh^2\theta}\right),} $
and consequently,
\begin{equation}
\left\{\begin{array}{l}
1\, \leq \, {\max_{\Sigma_1}\cosh\theta} \,
\leq \, 1+ c_{16}e^{-2c_{15}t}\\[2mm]
{\lambda_i^2(p,t)}\, \leq \, \frac{c_{16}e^{-2c_{15} t}}
{(1+ c_{16}e^{-2c_{15} t})}
\, \leq \, {c_{16}e^{-2c_{15} t}}=:{ (1-\delta(t))}
\end{array}\right.
\end{equation}
that is,  we  have for each $t$ a constant $\delta(t)$  explicitly defined,
and that approaches one  in an exponentially decreasing way,
and 
$$\frac d{dt}\cosh\theta \le \Delta \cosh\theta-\delta
(t)\cosh\theta ||B||^2.$$
 Setting $ p(t)=\frac{1}{\sqrt{nc_{16}}}e^{c_{15}t}$  and
 { $\psi= e^{\frac 12 c_{15}t}\cosh^{p(t)}\theta\|B\|^2,$} we have
\begin{eqnarray*}
\frac d{dt}\psi &\le& \Delta \psi-2\cosh ^{-p}\theta\, \nabla \cosh^p\theta
\,\nabla \psi 
-c_{17}\left\{e^{\frac 12 c_{15}t}\psi^2-e^{\frac
14c_{15}t}\psi^{\frac 12}-\psi\right\}.
\end{eqnarray*}
In \cite{lisa2} we show this
implies $\|B\|\leq c_{18}e^{-\tau t}$, where $\tau$ is a positive constant.
Since $F_t=F_0+\int_0^t H$ and the mean curvature is exponentially decreasing
we can conclude that $F_t(p)$ lies on a compact region of $\oM$,
and  for any sequence $t_N\to +\infty$
we obtain a subsequence $t_n$ such that  $F_{t_n}$ converges uniformly to
a spacelike graph of a map $f_{\infty}$. By (5) this map must be
constant. Furthermore, in this case  the
limit is the same, for any sequence $t_N\to +\infty$ we take.
This gives the next theorem, obtained with no need of using
Bernstein results:
\begin{theo}[{\cite{lisa2}}]
If $Ricci_1>0$ everywhere and $K_1\geq K_2$, $\Sigma_2$ not necessarily
compact,  
all the flow converges  to a unique slice.
\end{theo}
\section{Homotopy to a constant map}
We will give some applications of theorem 5.
We assume in this section $\Sigma_1$ is closed and $\Sigma_2$ is
complete with $R_2$ bounded and its derivatives. We also assume
either $K_1>0$ everywhere,  or  
$Ricci_1>0$ and  $K_2\leq -c< 0$ everywhere.

Given a constant $\rho>0$  we consider a new metric $g_1-g'_2$
on $\oM=\Sigma_1\times \Sigma_2$ where $g'_2=\rho^{-1}g_2$. Now
if $f:\Sigma_1\to \Sigma_2$ satisfies $f^*g_2<\rho g_1$,  means
$\Gamma_f$ is a timelike submanifold w.r.t. $g_1- g_2'$. Then
the curvature conditions in theorem 5 demands $K_1\geq \rho K_2$, that
 can be translated
in  the following 
\begin{theo}[{\cite{lisa2}}] There exist a constant 
$0\leq \rho\leq +\infty$,
such that any map $f:\Sigma_1\to \Sigma_2$ satisfying
$f^*g_2<\rho g_1$ is homotopically trivial.
If $K_1>0$ everywhere we may take $\rho\leq {\min_{\Sigma_1}K_1}/
{\sup_{\Sigma_2}K_2^+}$. For  $K_2\leq -c$ everywhere, we may take $\rho=+\infty$.
\end{theo}
\noindent
Note that, for $Ricci>0$, $K_2\leq -c<0$ everywhere,
then $\rho\geq {\max_{\Sigma_1}K_1^-}/{\inf_{\Sigma_2}-K_2}$,
where $K^-=\max\{-K, 0\}$. This means we may take $\rho=+\infty$ if
$K_2\leq -c$ as in case  $\sup_{\Sigma_2}K_2^+=0$ and $K_1>0$. 
This is the case $n=1$.
The homotopy is given by the flow, namely, since
$F_t(p)=(\phi_t(p),f_t(\phi_t(p)))$, where $\phi_t:\Sigma_1\to
\Sigma_1$ is a diffeomorphism with $\phi_0=id_{\Sigma_1}$, then
$K(t,p)=f_t(\phi_t(p))$ is the homotopy.
This gives a new proof of the classic Cartan-Hadmard theorem:
\begin{corol} If $K_2\leq 0$, $m\geq 2$, 
any map $f:\mathbb{S}^m\to \Sigma_2$
is homotopically trivial.
\end{corol}
\noindent
The condition given in \cite{w1},  $det(g_1+f^*g_2)<2$
implies
$\sum_i\lambda_i^2 + 1\leq \prod_i(1+\lambda_i^2) <2$
and so
$\Gamma_f$ is a spacelike
submanifold. The next theorem, obtained in the Riemannian
setting,  can be seen as a reformulated
corollary of theorem 5:
\begin{theo}[{\cite{w1,tw}}]  Assume both $\Sigma_i$ are closed
and with constant sectional curvature
$K_i$ and satisfying
$K_1\geq |K_2|$, $K_1+K_2>0$.\\
(1) If $det(g_1+f^*g_2)<2$, then $\Gamma_f$ can be deformed
 by a family of graphs to the one of a constant map.\\[1mm]
(2) If $f$ is an area decreasing map, that is
$\lambda_1\lambda_j<1$ for $i\neq j$, then it is homotopically trivial.
\end{theo}
The area decreasing condition is a slightly more general condition than
spacelike graph for $n\geq 2$. In case $n=1$ any map is area decreasing,
but it is included in the case $K_2\leq 0$.
We note that in the previous theorem 
it is used the Riemannian structure, and in this setting
$K_2$ cannot be given arbitrarily negative, a somehow artificial
condition, that can be dropped if one uses the pseudo-Riemannian
structure of the product.

\section*{Acknowledgements}
The first author is partially supported by NSFC (No.10501011) and by
Funda\c{c}\~{a}o Ci\^{e}ncia e Tecnologia (FCT) through a FCT
fellowship SFRH/BPD/26554/2006. The second author is
partially supported by FCT through the Plurianual of CFIF
and POCI-PPCDT/MAT/60671/2004.


\begin{thebibliography}{29}
\bibitem{aa} {A.\ Albujer A. and L.\ Al\'{i}as,
 \em Calabi-Bernstein results for
maximal  surfaces in Lorentz product spaces. \em 
arXiv:math/0709.4363, 2007 }
\bibitem{cal} {E.\ Calabi, \em Examples of Bernstein problems for some
nonlinear equations. \em Proc. Sympos. Pure Math. {\bf 15}(1970), 223--230.}
\bibitem{cy} {S.\ Cheng  and S.T.\ Yau, \em  Differential equations on
Riemannian manifolds and their geometric
 applications. \em  Comm.\ Pure Appl.\ Math.\ {\bf 28}(1975),  333--354.}
\bibitem{cheebin} { J.\ Cheeger and D.G.\ Ebin, \em Comparison theorems 
in Riemannian geometry. \em  North-Holland Mathematical Library, Vol. 9. 
North-Holland Publishing Co., Amsterdam-Oxford; 
American Elsevier Publishing Co., Inc., New York, 1975.}
\bibitem{chern} {S.S.\ Chern, \em Simple proofs of two theorems on 
minimal surfaces. \em  Enseignement Math. II. S\'{e}r {\bf 15} (1969), 53-61.}
\bibitem{eh} {K.\ Ecker and G.\ Huisken,
\em  Parabolic methods for the construction of spacelike slices of
prescribed mean curvature in cosmological spacetimes. \em  Comm.\ Math.\
Phys.\ {\bf 135} (1991), no. 3, 595--613.}
\bibitem{es}{ J.\ Eells and J.H.\ Sampson,
\em Harmonic mappings of Riemannian manifolds. 
\em Amer.\ J.\ Math.\ {\bf 86} (1964), 109--160..}
\bibitem{grom}{M.\ Gromov, \em Partial differential relations. \em
Ergebnisse der Mathematik und ihrer Grenzgebiete (3) 
[Results in Mathematics and Related Areas (3)], 9.
Springer-Verlag, Berlin, 1986.}
\bibitem{jx}{J.\ Jost and Y.\ Xin, \em  Some aspects of the global
geometry of entire spacelike submanifolds. \em  Results math. {\bf 40}(2001),
233--245.}
\bibitem{lisa1}{G.\ Li and I.M.C.\ Salavessa, \em  Graphic Bernstein
results in curved pseudo-Riemannian manifolds. \em  Arxiv.0801.3850. }
\bibitem{lisa2}{G.\ Li and I.M.C.\ Salavessa, \em  Mean curvature
flow of spacelike graphs. \em  Arxiv.0804.0783. }
\bibitem{sal}{I.M.C.\ Salavessa, \em Spacelike graphs with parallel 
mean curvature. \em  Bull.\ Bel.\ Math.\ Soc.\  {\bf 15 }
(2008), 65-76.}
\bibitem{tw}{M-P.\ Tsui, and M-T.\ Wang, \em 
Mean curvature flows and isotopy of
maps between spheres. \em 
 Comm.\ Pure Appl.\ Math.\  {\bf 57}  (2004),  no. 8, 1110--1126.}
\bibitem{w1}{M-T.\ Wang, 
\em Long-time existence and convergence of graphic
 mean curvature flow in arbitrary codimension. \em  Invent.\ Math.\ 
{\bf 148} (2002), no. 3, 525-543.}
\bibitem{xin}{Y.\ Xin, \em Mean curvature flow with convex
Gauss image. \em   Chin.\ Ann.\ Math.\ Ser.\ B {\bf 29} (2008), no. 2, 
121--134.}
\bibitem{z}{X-P.\ Zhu, \em Lectures on mean curvature flows. \em 
AMS/IP Studies in advanced mathematics, {\bf 32}, 
American Mathematical Society, Providence, RI; 
International Press (2002).}
\end{thebibliography}
\end{document}